# STATIONARITY AND GEOMETRIC ERGODICITY OF A CLASS OF NONLINEAR ARCH MODELS

By Youssef Saïdi and Jean-Michel Zakoïan

*Université Mohammed V, Université Lille 3 and CREST*

A class of nonlinear ARCH processes is introduced and studied. The existence of a strictly stationary and $\beta$-mixing solution is established under a mild assumption on the density of the underlying independent process. We give sufficient conditions for the existence of moments. The analysis relies on Markov chain theory. The model generalizes some important features of standard ARCH models and is amenable to further analysis.

**1. Introduction.** Since the appearance of seminal papers by Engle [9] and Bollerslv [2], a variety of GARCH (generalized autoregressive conditionally heteroskedastic) specifications have been introduced to model the characteristic features of observed financial time series. These specifications are of the form

$$\varepsilon_t = \sigma_t \eta_t, \qquad t \in \mathbb{Z}, \tag{1.1}$$

where the sequence $(\eta_t)$ is independent and identically distributed (i.i.d.) with zero mean and unit variance, and $\sigma_t$ is a positive variable called *volatility*, which is a measurable function of the past, $\{\varepsilon_{t-i}, i > 0\}$. Typically, $\varepsilon_t$ represents the logarithm of the return, that is, the variation of the price in logarithm.

The original model specified $\sigma_t^2$ as a linear function of the squared past log-returns and was found adequate to capture many *stylized facts* associated with the financial data, namely tail heaviness, volatility clustering, leptokurtosis of the marginal distribution and dependence without autocorrelation. Other characteristic properties such as asymmetries motivated extensions of the basic model (see, e.g., [15, 20]). A common feature of these models is that $\sigma_t$ is specified as a strictly increasing function of the modulus of the past returns. In general, the specification of $\sigma_t$ involves a linear combination of some function of the past returns.

In this paper, we consider a class of *nonlinear* ARCH processes. More precisely, the model we study in this paper is given by

$$
\begin{aligned}
\varepsilon_t &= \sigma_t \eta_t, \\
\sigma_t^2 &= \omega + \alpha \varepsilon_{t-1}^2 \mathbb{1}_{\varepsilon_{t-1}^2 > k \varepsilon_{t-2}^2},
\end{aligned}
\tag{1.2}
$$

where $\omega$, $\alpha$ and $k$ are nonnegative constants with $\omega > 0$ and where the same assumptions are made concerning $(\eta_t)$ as in (1.1). The standard ARCH(1) model is obtained as a particular case by taking $k = 0$. The conditionally homoskedastic model (constant volatility) can be obtained by setting $\alpha = 0$, but it is worthnoting that "large" values of $k$ also produce a model which is close to being homoskedastic. This model belongs to the class of endogenous switching regime models, in the spirit of the threshold autoregressive models of Tong and Lim [19]. In the present model, the volatility equation can be interpreted as a two-regime specification, the first regime being homoskedastic ($\sigma_t^2 = \omega$) and the second one being a classical ARCH(1) ($\sigma_t^2 = \omega + \alpha \varepsilon_{t-1}^2$). The originality of the specification, however, is that the regime change depends on the *relative* variation of the last squared observation. As soon as the relative variations ($\varepsilon_{t-1}^2/\varepsilon_{t-2}^2$) are small, the process remains in the homoskedastic regime. But, when these variations are large, the volatility depends on the last squared observation. The coefficient $k$ allows for flexibility in the occurrence of the two regimes. Empirical motivations for model (1.2) based on the features of real financial time series can be found in the dissertation by Saïdi [17].

The aim of this paper is to study the stability properties of the specification in (1.2). Recent references dealing with ergodic properties of GARCH-type models are [1, 5, 10, 11]. These papers use a random coefficient linear representation of the volatility, of the form $\sigma_t^2 = \omega(\eta_{t-1}) + a(\eta_{t-1})\sigma_{t-1}^2$ in the first-order case, which does not hold in our framework. A different approach is used by Cline and Pu [7] who establish sharp conditions for geometric ergodicity of a class of threshold autoregressive ARCH models under assumptions we will discuss further.

The rest of the paper proceeds as follows. In Section 2, we recall the main results of Markov chain theory that we will use in the sequel. Section 3 is devoted to the existence of strictly stationary solutions. We start by considering the deterministic model implied by (1.2). Then we establish conditions for the existence of strictly stationary and $\beta$-mixing solutions. Finally, we provide conditions for the existence of moments.

**2. Some Markov chain results.** In this section, we give results from the theory of Markov chain processes that allow to study the existence of ergodic solutions to stochastic difference equations. This section is heavily based on the book by Meyn and Tweedie [13]. Let $E \subset \mathbb{R}^d$ and let $\mathcal{E}$ be the Borel $\sigma$-field on $E$. We denote by $\{X_t, t \geq 0\}$ a homogeneous Markov chain on $(E, \mathcal{E})$ and denote by $P^t(x, B) = \mathbb{P}(X_t \in B | X_0 = x)$ the probability of moving from $x \in E$ to the set $B \in \mathcal{E}$ in $t$ steps. The Markov chain $(X_t)$ is $\phi$-*irreducible* if, for some nontrivial $\sigma$-finite measure $\phi$ on $(E, \mathcal{E})$,

$$\forall B \in \mathcal{E} \quad \phi(B) > 0 \implies \forall x \in E, \exists t > 0, \quad P^t(x, B) > 0.$$

If $(X_t)$ is $\phi$-irreducible, there exists a *maximal irreducibility measure $M$* (see [13], Proposition 4.2.2) and we set $\mathcal{E}^+ = \{B \in \mathcal{E} | M(B) > 0\}$. The chain is called *positive recurrent* if

$$\forall x \in \mathcal{E}, \ \forall B \in \mathcal{E}^+ \quad \limsup_{t \to \infty} P^t(x, B) > 0.$$

For a $\phi$-irreducible Markov chain, positive recurrence is equivalent (see [13], Theorem 18.2.2) to the existence of a (unique) *invariant distribution*, that is, a probability measure $\pi$ such that

$$\forall B \in \mathcal{B} \quad \pi(B) = \int P(x, B) \pi(dx).$$

"Geometric ergodicity" refers to the rate of convergence of the transition probabilities to the invariant distribution. More precisely, if $\|\cdot\|$ denotes the total variation norm, the Markov chain $(X_t)$ is said to be *geometrically ergodic* if there exists a $\rho$, $\rho \in (0,1)$ such that

(2.1) $\quad \forall x \in E \quad \rho^{-t} \|P^t(x, \cdot) - \pi\| \to 0 \quad \text{as } t \to +\infty.$

In order to state the following criterion for the geometric ergodicity of a Markov chain, we need the notions of $T$-chain, small sets and aperiodicity. For any distribution $a = (a_n)$ on the set of positive integers, for all $x \in E$ and $B \in \mathcal{E}$, let $K_a(x, B) = \sum_{n \geq 1} a_n P^n(x, B)$. Recall that if $E$ is endowed with a metric, a function $h: E \to \mathbb{R}$ is called *lower semicontinuous* if for any constant $c$, the set $\{x : h(x) > c\}$ is open. Now, if for any open set $B$, the function $P(\cdot, B)$ is lower semicontinuous, $(X_t)$ is called a *Feller* Markov chain. More generally, if there exists a function $T: E \times \mathcal{E} \to [0, +\infty)$ and a distribution $a = (a_n)$ on the set of positive integers such that (i) $T(\cdot, B)$ is lower semicontinuous, $\forall B \in \mathcal{E}$, (ii) $T(x, \cdot)$ is a nontrivial measure over $(E, \mathcal{E})$, $\forall x \in E$ and (iii) $K_a(x, B) \geq T(x, B), \forall x \in E, B \in \mathcal{E}$, then $(X_t)$ is called a *T-chain* and $T$ is called a *continuous component* of $K_a$. A set $C \in \mathcal{E}$ is called a $\nu_m$-*small set* if there exist an $m > 0$ and a nontrivial measure $\nu_m$ on $\mathcal{E}$ such that $\forall x \in C$ and $\forall B \in \mathcal{E}, P^m(x, B) \geq \nu_m(B)$. Let $C$ be a $\nu_M$-small set where the measure $\nu_M := \nu$ is such that $\nu(C) > 0$. Such a measure exists whenever $C \in \mathcal{E}^+$ (see [13], Proposition 5.2.4). Let $E_C = \{m \geq 1 | C$ is $\nu_m$-small with $\nu_m = \delta_m \nu$ for some $\delta_m > 0\}$. Then if $(X_t)$ is a $\phi$-irreducible Markov chain and $C \in \mathcal{B}^+$, the greatest common divisor $d$ of the set $E_C$ does not depend on $C$ and is called the *period* of the Markov chain. If $d = 1$, $(X_t)$ is said to be *aperiodic*. If every compact set is small, then $(X_t)$ is a $T$-chain. If $(X_t)$ is a $\phi$-irreducible $T$-chain, then every compact set is small (see [13], Proposition 5.5.7 and Theorem 6.2.5). However, some noncompact sets may also be small, and such sets can be worth considering, as we shall see.

We are now in a position to state a criterion for geometric ergodicity based on $m$-step transitions, which is adapted from [13], Theorem 19.1.3. The use of $m$-step transitions in ergodicity criteria was suggested by Tjøstheim [18].

THEOREM 2.1. *Assume that:*

  (i) $(X_t)$ *is $\phi$-irreducible for some measure $\phi$ on $(E, \mathcal{E})$,*
 (ii) $(X_t)$ *is an aperiodic $T$-chain,*
(iii) *there exists a small set $C \in \mathcal{E}^+$, an integer $m \geq 1$ and a nonnegative continuous function (test function) $g: E \to [0, +\infty)$ such that*

$$E[g(X_{t+m})|X_t = x] \leq \begin{cases} (1-\beta)g(x) - \beta, & x \in C^c, \\ b, & x \in C, \end{cases}$$

*for some strictly positive constants $\beta$ and $b$. Then $(X_t)$ is geometrically ergodic. Moreover, $E_\pi g(X_t)$ is finite, where $E_\pi$ denotes expectation taken under the stationary distribution.*

One consequence of the geometric ergodicity is that the Markov chain $(X_t)$ is $\beta$-*mixing*, and hence *strongly mixing, with geometric rate*. Recall that for a stationary process, the $\beta$-mixing coefficients are defined by

$$(2.2) \qquad \beta_X(k) = E \sup_{B \in \sigma(X_s, s \geq k)} |\mathbb{P}(B|\sigma(X_s, s \leq 0)) - \mathbb{P}(B)|.$$

The process is called $\beta$-mixing if $\lim_{k \to \infty} \beta_X(k) = 0$. If $Y = (Y_t)$ is a process such that $Y_t = f(X_t, \ldots, X_{t-r})$ for some measurable function $f$ and some integer $r \geq 0$, then $\sigma(Y_t, t \leq s) \subset \sigma(X_t, t \leq s)$ and $\sigma(Y_t, t \geq s) \subset \sigma(X_{t-r}, t \geq s)$. Thus,

$$(2.3) \qquad \beta_Y(k) \leq \beta_X(k - r) \qquad \text{for all } k \geq r.$$

Davydov [8] showed that for an ergodic Markov chain $(X_t)$ with invariant probability measure $\pi$,

$$\beta_X(k) = \int \|P^k(x, \cdot) - \pi\| \pi(dx).$$

Noting that in (2.1), the rate $\rho$ can be chosen independently of the initial point $x$, it follows that $\beta_X(k) = O(\rho^k)$ if (2.1) holds.

**3. Existence of stationary solutions.** In this section, we consider the problem of the existence of strictly stationary and second-order stationary solutions to model (1.2). The problem is not standard because, contrary to most ARCH-type specifications, no linear representation of the model seems to exist. Hence, we cannot rely on the theory developed in the papers by Bougerol and Picard [3, 4]. Instead, we will use the techniques of Tweedie to deal with the stationarity question.

Thinking of the standard ARCH(1) model, we could perhaps expect to require a (strict and second-order) stationarity condition of the form $\alpha < 1$. The presence of the (conditionally) homoskedastic regime seems to allow us greater freedom. As for the threshold autoregressive models, it will be helpful to first consider the deterministic model.

3.1. *Stability of the deterministic model.* Suppose that, in model (1.2), the i.i.d. process $(\eta_t)$ is such that $\eta_t^2 = 1$, for all $t$, almost surely. We call this model deterministic, although the sign of $\varepsilon_t$ is, of course, a random variable. For ease of exposition, we take $\varepsilon_0 = 0$, but any other initial value would also produce the following asymptotic results:

THEOREM 3.1. *Let $(\varepsilon_t)_{t \geq 0}$ be as defined in (1.2), with $\eta_t^2 = 1$, a.s. for all $t$, and $\varepsilon_0 = 0$. Then:*

(i) *if $\max(\alpha, 1) < k$ or $\alpha = 0$, then there exists $i \geq 3$ such that $\forall t \geq i$, $\varepsilon_t^2 = \omega$ a.s.;*

(ii) *if $\alpha < 1$ and $k \leq 1$, then $\varepsilon_t^2 \longrightarrow \frac{\omega}{1-\alpha}$ a.s. when $t \to +\infty$;*

(iii) *if $\alpha \geq \max(1, k)$, then $\varepsilon_t^2 \longrightarrow +\infty$ a.s. when $t \to +\infty$.*

PROOF. We have, a.s., $\varepsilon_1^2 = \omega$ and $\varepsilon_2^2 = \omega(1+\alpha)$. The value of $\varepsilon_3^2$ depends on the position of $1 + \alpha$ compared to $k$. Let, for all $i \geq 0$,

$$E_i = \left\{ (\alpha, k); \alpha + \frac{1}{1 + \alpha + \cdots + \alpha^i} \leq k \right\}.$$

Since $\alpha \geq 0$, the sets $E_i$ constitute an increasing sequence. We have

$$E_\infty := \bigcup_{i \geq 0} E_i = \{\alpha = 0, k \geq 1\} \cup \{0 < \alpha < 1 < k\} \cup \{1 \leq \alpha < k\}$$

$$= \{\max(\alpha, 1) < k\} \cup \{(0, 1)\}.$$

Let us consider the different cases.

*Case* (i). We have $(\alpha, k) \in E_\infty$, hence there exists $i \geq 0$ such that $(\alpha, k) \in E_i$. Let $i_0 = \min\{i \geq 0, (\alpha, k) \in E_i\}$. For $1 \leq i \leq i_0 + 2$, we have $\varepsilon_i^2 = \omega(1 + \cdots + \alpha^{i-1})$. Then

$$\frac{\varepsilon_{i_0+2}^2}{\varepsilon_{i_0+1}^2} = \alpha + \frac{1}{1 + \alpha + \cdots + \alpha^{i_0}} \leq k.$$

It follows that $\varepsilon_{i_0+3}^2 = \omega$ and, since $\frac{\varepsilon_{i_0+3}^2}{\varepsilon_{i_0+2}^2} = \frac{1}{1+\alpha+\cdots+\alpha^{i_0+1}} \leq 1 < k$, that $\varepsilon_i^2 = \omega$ for all $i \geq i_0 + 3$.

*Case* (ii). If $(\alpha, k) \neq (0, 1)$, then $(\alpha, k) \notin E_\infty$. Thus, for all $i \geq 1$, $\varepsilon_i^2 = \omega(1 + \cdots + \alpha^{i-1})$ and the result follows. When $(\alpha, k) = (0, 1)$, the sequence $(\varepsilon_i^2)$ takes the constant value $\omega$.

*Case* (iii). We have $(\alpha, k) \notin E_\infty$. Thus, for all $i \geq 1$, $\varepsilon_i^2 = \omega(1 + \cdots + \alpha^{i-1})$ and the sequence $(\varepsilon_i^2)$ tends to $+\infty$. □

From this result, the region of nonexplosion of the deterministic models is given by $\alpha < \max(1, k)$. We now turn to the general case.

3.2. *Markov chain results.* As with many discrete-time models, the analysis of the probability structure of model (1.2) draws on Markov chain results. Let

$$X_t = \begin{pmatrix} \varepsilon_t^2 \\ \varepsilon_{t-1}^2 \end{pmatrix} = \begin{pmatrix} X_{1,t} \\ X_{2,t} \end{pmatrix}$$

and let

$$\forall x \in \mathbb{R}^2 \qquad \psi(x) = \omega + \alpha x_1 \mathbb{1}_{x_1 > kx_2}.$$

The vector representation of model (1.2) takes the form of a nonlinear stochastic difference equation,

$$(3.1) \qquad X_t = \begin{pmatrix} \psi(X_{t-1})\eta_t^2 \\ X_{1,t-1} \end{pmatrix} := F(X_{t-1}, \eta_t), \qquad t \geq 1,$$

where the i.i.d. sequence $(\eta_t)$ is supposed to be independent of the initial state $X_0$. Note that models of the form (3.1) are considered, among others, by [13], Chapter 7, but under a smoothness assumption on the function $F$ which is not valid in our framework. Let $\lambda_m^+$ be the Lebesgue measure and let $\mathcal{B}(\mathbb{R}^{+m})$ be the Borel class of sets for $\mathbb{R}^{+m}$. We will make the following assumption:

ASSUMPTION A. *The variables $\eta_t^2$ admit a density $f$ with respect to $\lambda_1^+$, with $f > 0$ on $\mathbb{R}^+$. Moreover, $E\eta_t = 0$ and $E\eta_t^2 = 1$.*

LEMMA 1. *The process $(X_t)_{t \geq 0}$ is a time-homogeneous Markov chain on $\mathbb{R}^{+2}$, with transition probabilities given as follows:*

$$\forall x = (x_1, x_2) \in \mathbb{R}^{+2}, \forall B = B_1 \times B_2 \in \mathcal{B}(\mathbb{R}^{+2}).$$

$$(3.2) \qquad \mathbf{P}(x, B) = P[\eta_t^2 \in \psi(x)^{-1} B_1] \mathbb{1}_{x_1 \in B_2}.$$

*Moreover, under Assumption A, the process $(X_t)$ is $\lambda_2^+$-irreducible ($\lambda_2^+$ is therefore a maximal irreducibility measure).*

PROOF. Equation (3.1) ensures that $(X_t)$ is a time-homogeneous Markov chain. The two-step transition probabilities are given as follows:

$$\forall B = B_1 \times B_2 \in \mathcal{B}(\mathbb{R}^{+2}), \forall x \in \mathbb{R}^2$$

$$\mathbf{P}^2(x, B) = P[\psi(X_{t-1})\eta_t^2 \in B_1, \psi(x)\eta_{t-1}^2 \in B_2 | X_{t-2} = x]$$

(3.3)
$$= P[\psi\{\psi(x)\eta_{t-1}^2, x_1\}\eta_t^2 \in B_1, \psi(x)\eta_{t-1}^2 \in B_2]$$

$$= \int \mathbb{1}_{\psi(x)^{-1}B_2}(y) P[\eta_t^2 \in \psi\{\psi(x)y, x_1\}^{-1}B_1] f(y)\, d\lambda_1^+(y).$$

This can be seen by using the Fubini theorem, using the independence between $\eta_t$ and $\eta_{t-1}$ and noting that $\psi(\cdot) \geq \omega > 0$. If $\lambda_1^+(B_1) > 0$, we have, for all $y \in \mathbb{R}^+$, $P[\eta_t^2 \in \psi\{\psi(x)y, x_1\}^{-1}B_1] > 0$, in view of Assumption A. Similarly, $\lambda_1^+\{\psi(x)^{-1}B_2\} > 0$ if $\lambda_1^+(B_2) > 0$. Hence, $\mathbf{P}^2(x, B) > 0$, which ensures that $\lambda_2^+$ is an irreducibility measure. We have $\mathbf{P}^t(x, B) = 0$ for any Borel set $B \subset (\mathbb{R}^-)^2$, any $t > 0$ and any $x \in \mathbb{R}^2$. Thus, any irreducibility measure $\phi$ is such that $\phi(B) = 0$ for any $B \in \mathcal{B}(\mathbb{R}^{-2})$. It follows that $\lambda_2^+$ is a maximal irreducibility measure (see [13], Proposition 4.2.2). □

REMARK 1. Cline and Pu [6] provide conditions for irreducibility (as well as aperiodicity and the $T$-chain property) for a general class of nonlinear autoregressive models encompassing (1.2). Since we use slightly weaker conditions for the error density, we give direct proofs of the corresponding lemmas.

REMARK 2. The transition probability defined in (3.2) is a function of $x$ which is not lower semicontinuous for any open set $B$. To see this, let $x_1 = kx_2$, let $B = B_1 \times B_2$ be an open set such that $p = P[\eta_t^2 \in \omega^{-1}B_1] > P[\eta_t^2 \in (\omega + \alpha x_1)^{-1}B_1] = q$ and such that $x_1 \in B_2$. For $x = (x_1, x_2)$, we have $P[\eta_t^2 \in \psi(x)^{-1}B_1] > c = (p+q)/2$. Any neighborhood of $x$ contains points $y = (y_1, y_2)$ with $y_1 > ky_2$. For such points, we have $\psi(y) = \omega + \alpha y_1$ and thus, if $y_1$ is sufficiently close to $x_1$, $P[\eta_t^2 \in \psi(y)^{-1}B_1] < c$. The set $\{x : \mathbf{P}(x, B) > c\}$ is therefore not open. It follows that $(X_t)$ is not a Feller chain. The fact that compact sets are small, which will be used in the verification of our ergodicity criterion, is thus not straightforward. This property will follow from the next result.

LEMMA 2. *Under the assumptions of Lemma 1, the process $(X_t)$ is a $T$-chain.*

PROOF. It will be convenient to consider a partition of the positive quadrant of $\mathbb{R}^2$ into three regions: $D_1 = \{x_1 < kx_2\}$, $D_2 = \{x_1 = kx_2\}$ and $D_3 = \{x_1 > kx_2\}$.

For $x \in D_1 \cup D_2$, we have $\psi(x) = \omega$. Thus, from (3.3), using the Fubini theorem and the independence between $\eta_t$ and $\eta_{t-1}$, we have

$$\forall x \in D_1 \cup D_2, \forall B = B_1 \times B_2 \in \mathcal{B}(\mathbb{R}^{+2})$$

$$\mathbf{P}^2(x, B) = P[\psi(\omega \eta_{t-1}^2, x_1) \eta_t^2 \in B_1, \; \omega \eta_{t-1}^2 \in B_2]$$

$$= \int \mathbb{1}_{\omega^{-1} B_2 \cap (-\infty, \omega^{-1} k x_1]}(y) P[\eta_t^2 \in \omega^{-1} B_1] f(y) \, d\lambda_1^+(y)$$

$$+ \int \mathbb{1}_{\omega^{-1} B_2 \cap (\omega^{-1} k x_1, +\infty)}(y) P[\eta_t^2 \in \{\omega(1 + \alpha y)\}^{-1} B_1] f(y) \, d\lambda_1^+(y).$$

By the Lebesgue theorem and Assumption A, we can conclude that $\mathbf{P}^2(\cdot, B)$ is continuous over the set $D_1$. This is not the case for $x \in D_2$. However, if some sequence $(x_n)$ converges to $x$ with $x_n \in D_1 \cup D_2$, we have $\mathbf{P}^2(x_n, B) \longrightarrow \mathbf{P}^2(x, B)$ by the same arguments. For $x_n = (x_{1n}, x_{2n}) \in D_3$, we have

$$\mathbf{P}^2(x_n, B) = P[\psi\{(\omega + \alpha x_{1n})\eta_{t-1}^2, x_{1n}\}\eta_t^2 \in B_1, (\omega + \alpha x_{1n})\eta_{t-1}^2 \in B_2],$$

because $\psi(x_n) = \omega + \alpha x_{1n}$. Therefore, proceeding as for $D_1$,

$$\lim_{x_n \to x, x_n \in D_3} \mathbf{P}^2(x_n, B) = P[\psi\{(\omega + \alpha x_1)\eta_{t-1}^2, x_1\}\eta_t^2 \in B_1,$$

$$(\omega + \alpha x_1)\eta_{t-1}^2 \in B_2].$$

Setting

$$T(x, B) = P[\psi\{\omega \eta_{t-1}^2, x_1\}\eta_t^2 \in B_1, \omega \eta_{t-1}^2 \in B_2,$$

$$\psi\{(\omega + \alpha x_1)\eta_{t-1}^2, x_1\}\eta_t^2 \in B_1, (\omega + \alpha x_1)\eta_{t-1}^2 \in B_2],$$

we define a measure for any $x$, which is nontrivial because $T(x, \mathbb{R}^{+2}) = 1$. Setting $a(x_1) = \omega + \alpha x_1$, $T(x, B)$ can be decomposed into three probabilities, depending on the position of $\eta_{t-1}^2$, as follows:

$$P[\omega \eta_t^2 \in B_1, \omega \eta_{t-1}^2 \in B_2, a(x_1) \eta_{t-1}^2 \in B_2, \eta_{t-1}^2 < k x_1 / a(x_1)]$$

$$+ P[\omega \eta_t^2 \in B_1, \omega \eta_{t-1}^2 \in B_2, \{\omega + a(x_1) \eta_{t-1}^2\} \eta_t^2 \in B_1,$$

$$a(x_1) \eta_{t-1}^2 \in B_2, \eta_{t-1}^2 \in [k x_1 / a(x_1), k x_1 / \omega]]$$

$$+ P[\omega(1 + \alpha \eta_{t-1}^2) \eta_t^2 \in B_1, \omega \eta_{t-1}^2 \in B_2,$$

$$\{\omega + a(x_1) \eta_{t-1}^2\} \eta_t^2 \in B_1, a(x_1) \eta_{t-1}^2 \in B_2, \eta_{t-1}^2 > k x_1 / \omega].$$

This, in view of Assumption A, shows that the function $T(\cdot, B)$ is continuous. Finally, $\mathbf{P}^2(x, B) \geq T(x, B)$ for all $x$ and all $B$. Thus, $T$ is a continuous component of $\mathbf{P}^2$. The conclusion follows. □

Classical ergodicity proofs for nonlinear stochastic difference equations (as, for instance, in the case of TAR models, see [16]) rely on verifying a drift condition when the chain goes outside a compact set. In the model of this paper, no drift condition holds over the region $\{x_1 \leq kx_2\}$. It is therefore necessary to consider more general small sets than compact sets, as was done, for instance, by Cline and Pu [6], Theorem 2.5.

LEMMA 3. *Under the assumptions of Lemma 1, the set $C = \{x_1 \leq kx_2\}$ is small for the Markov chain $(X_t)$. Moreover, the chain is aperiodic.*

PROOF. For $x = (x_1, x_2) \in C$, we have $\psi(x) = \omega$. Thus, by (3.3), for any $B = B_1 \times B_2 \in \mathcal{B}(\mathbb{R}^{+2})$

$$
\begin{aligned}
\mathbf{P}^2(x, B) &= P[\psi(\omega \eta_{t-1}^2, x_1)\eta_t^2 \in B_1, \omega \eta_{t-1}^2 \in B_2] \\
&= P[\omega \eta_t^2 \in B_1, \omega \eta_{t-1}^2 \leq kx_1, \omega \eta_{t-1}^2 \in B_2] \\
&\quad + P[(\omega + \alpha \omega \eta_{t-1}^2)\eta_t^2 \in B_1, \omega \eta_{t-1}^2 > kx_1, \omega \eta_{t-1}^2 \in B_2] \\
&:= \mathbf{P}_1(x, B) + \mathbf{P}_2(x, B).
\end{aligned}
$$
(3.4)

Let $\varepsilon > 0$. For $x_1 > \varepsilon$, we have

(3.5) $\quad \mathbf{P}_1(x, B) \geq P[\omega \eta_t^2 \in B_1, \omega \eta_{t-1}^2 \leq k\varepsilon, \omega \eta_{t-1}^2 \in B_2] := \mu_1(B).$

For $x_1 \leq \varepsilon$, we have

$$\mathbf{P}_2(x, B) \geq P[(\omega + \alpha \omega \eta_{t-1}^2)\eta_t^2 \in B_1, \omega \eta_{t-1}^2 > k\varepsilon, \omega \eta_{t-1}^2 \in B_2] := \mu_2(B).$$

The measures $\mu_1$ and $\mu_2$ are clearly nontrivial. It follows that the sets $C_1 = \{x_1 > \varepsilon\} \cap C$ and $C_2 = \{x_1 \leq \varepsilon\} \cap C$ are small. The union of two small sets being a small set, we may conclude that $C = C_1 \cup C_2$ is a small set.

To prove aperiodicity, we will consider three-step transition probabilities. Recall that, for a $\phi$-irreducible Markov chain, the definition of the period $d$ is independent of the choice of a small set. For our small set, we choose $C_1$. For $x \in C_1$ and for $B = B_1 \times B_2 \in \mathcal{B}(\mathbb{R}^{+2})$, we have, from (3.4) and (3.5), after translation of the times,

$$
\begin{aligned}
\mathbf{P}^2(x, B) &\geq P[\omega \eta_{t+1}^2 \in B_1, \omega \eta_t^2 \leq k\varepsilon, \omega \eta_t^2 \in B_2] \\
&\geq P[\omega \eta_{t+1}^2 \in B_1, \omega \eta_{t-1}^2 \leq k\varepsilon, \omega \eta_t^2 \leq k\varepsilon, \omega \eta_t^2 \in B_2, \eta_t^2 \leq k\eta_{t-1}^2] \\
&:= \mu(B), \\
\mathbf{P}^3(x, B) &= P[\psi\{\psi(\omega \eta_{t-1}^2, x_1)\eta_t^2, \omega \eta_{t-1}^2\}\eta_{t+1}^2 \in B_1, \psi(\omega \eta_{t-1}^2, x_1)\eta_t^2 \in B_2] \\
&= P[\psi\{\omega \eta_t^2, \omega \eta_{t-1}^2\}\eta_{t+1}^2 \in B_1, \omega \eta_{t-1}^2 \leq kx_1, \omega \eta_t^2 \in B_2] \\
&\quad + P[\psi\{(\omega + \alpha \omega \eta_{t-1}^2)\eta_t^2, \omega \eta_{t-1}^2\}\eta_{t+1}^2 \in B_1, \\
&\qquad \omega \eta_{t-1}^2 > kx_1, (\omega + \alpha \omega \eta_{t-1}^2)\eta_t^2 \in B_2] \\
&\geq P[\psi\{\omega \eta_t^2, \omega \eta_{t-1}^2\}\eta_{t+1}^2 \in B_1, \omega \eta_{t-1}^2 \leq k\varepsilon, \omega \eta_t^2 \in B_2]
\end{aligned}
$$

$$\geq P[\omega\eta_{t+1}^2 \in B_1, \omega\eta_{t-1}^2 \leq k\varepsilon, \omega\eta_t^2 \in B_2, \eta_t^2 \leq k\eta_{t-1}^2]$$
$$\geq \mu(B).$$

The set $C_1$ is then both $\nu_2$-small and $\nu_3$-small, where $\nu_2 = \nu_3 = \mu$. This measure $\mu$ is nontrivial. The greatest common divisor $d$ of the set $E_{C_1}$ which appears in the definition of periodicity is thus equal to 1. The conclusion follows. $\square$

3.3. *$\beta$-mixing.* The main result of this paper is the following theorem:

THEOREM 3.2. *Under Assumption* A *and the condition* $k > 0$, *there exists a strictly stationary solution* $(\varepsilon_t)$ *to model* (1.2). *This solution is $\beta$-mixing, and hence strongly mixing, with geometric rate. Moreover, there exists* $r > 0$ *such that* $E_\pi(\varepsilon_t^{2r}) < \infty$.

REMARK 3. It is worth noting that when $k > 0$, strict stationarity holds regardless of the value of $\alpha$. When $k = 0$, that is, in the case of the standard ARCH(1), we have the well-known strict stationarity condition established by Nelson [14]: $0 \leq \alpha < \exp\{-E(\log \eta_t^2)\}$.

REMARK 4. Assumption A is crucial for strict stationarity to hold without an upper bound for $\alpha$. For instance, in the deterministic case, $\eta_t^2 = 1$, a.s., Assumption A is not verified and it was seen in Section 3.1 that stability requires $k > \alpha$, or $k \leq \alpha < 1$.

REMARK 5. Cline and Pu [7] provided useful conditions for geometric ergodicity of a general class of nonlinear AR–ARCH models. We cannot rely on their results, however, because in particular their Assumption A.5 does not hold for model (1.2).

To prove Theorem 3.2, we start by establishing the following lemma:

LEMMA 4. *Under the assumptions of Theorem* 3.2, *the Markov chain* $(X_t)$ *is geometrically ergodic.*

PROOF. The conclusion being obvious when $\alpha = 0$, we consider the case $\alpha > 0$. The proof consists in verifying the three conditions of Theorem 2.1 for $m = 2$. Property (i) holds with $\phi = \lambda_2^+$, by Lemma 1, (ii) holds by Lemmas 2 and 3. To check (iii), we take $g(x) = g(x_1, x_2) = x_1^r$, where $r \in (0, 1]$. Let $\mu_{2r} = E(\eta_t^{2r})$ and let $\mu_{2r}^* = E(\eta_t^{2r} \mathbb{1}_{\eta_t^2 > k/\alpha})$. Note that these quantities are finite under Assumption A. We have

$$
\begin{aligned}
E[g(X_{t+2})|X_t &= (x_1, x_2)] \\
&= E[\varepsilon_{t+2}^{2r}|X_t = (x_1, x_2)] \\
&= E[\eta_{t+2}^{2r}(\psi\{\psi(x_1, x_2)\eta_{t+1}^2, x_1\})^r] \\
&= \mu_{2r} E[\psi\{\psi(x_1, x_2)\eta_{t+1}^2, x_1\}]^r \\
&= \mu_{2r} E\{\omega + \alpha\psi(x_1, x_2)\eta_{t+1}^2 \mathbb{1}_{\eta_{t+1}^2 > kx_1/\psi(x_1,x_2)}\}^r \\
&\leq \mu_{2r}\omega^r + \mu_{2r}\alpha^r \psi(x_1, x_2)^r E[\eta_{t+1}^{2r} \mathbb{1}_{\eta_{t+1}^2 > kx_1/\psi(x_1,x_2)}],
\end{aligned}
\tag{3.6}
$$

where the last inequality follows from the elementary inequality $(a+b)^r \leq a^r + b^r$ for any $a, b \geq 0$. For $x_1 > kx_2$, we then have

$$E[g(X_{t+2})|X_t = (x_1, x_2)] \leq \mu_{2r}\omega^r + \mu_{2r}\alpha^r(\omega + \alpha x_1)^r E[\eta_{t+1}^{2r} \mathbb{1}_{\eta_{t+1}^2 > kx_1/(\omega+\alpha x_1)}].$$

When $x_1 \to +\infty$, the right-hand side of this inequality is equivalent to

$$\alpha^{2r}\mu_{2r}\mu_{2r}^* x_1^r.$$

Now $\alpha^{2r}\mu_{2r}\mu_{2r}^*$ tends to $P[\eta_t^2 > k/\alpha]$ when $r \to 0$, by the Lebesgue theorem. This probability being strictly less than 1 when $k > 0$ (by Assumption A), we have $\alpha^{2r}\mu_{2r}\mu_{2r}^* < 1$ for $r$ sufficiently small. Therefore, there exist $\beta > 0$, $r > 0$ and $M > 0$ such that

$$x_1 > M \quad \text{and} \quad x_1 > kx_2 \quad \Longrightarrow \quad E[\varepsilon_{t+2}^{2r}|X_t = (x_1, x_2)] \leq (1-\beta)x_1^r - \beta.$$

For $x_1 \leq M$, we have $\psi(x_1, x_2) \leq \omega + \alpha M$ and hence, from

$$E[\eta_{t+1}^{2r}\mathbb{1}_{\eta_{t+1}^2 > kx_1/\psi(x_1,x_2)}] \leq \mu_{2r}$$

and (3.6), we have

$$E[\varepsilon_{t+2}^{2r}|X_t = (x_1, x_2)] \leq \mu_{2r}\omega^r + \mu_{2r}^2\alpha^r(\omega + \alpha M)^r.$$

Finally, for $x_1 \leq kx_2$, since $\psi(x_1, x_2) = \omega$, we have, by (3.6),

$$E[\varepsilon_{t+2}^{2r}|X_t = (x_1, x_2)] \leq \mu_{2r}\omega^r(1 + \mu_{2r}\alpha^r) \leq \mu_{2r}\omega^r + \mu_{2r}^2\alpha^r(\omega + \alpha M)^r.$$

We can conclude that (iii) holds, with $C = [0, M]^2 \cup \{x_1 \leq kx_2\}$ and $b = \mu_{2r}\omega^r + \mu_{2r}^2\alpha^r(\omega + \alpha M)^r$.

That $C$ is a small set is a consequence of Lemma 2 (implying that any compact set is small), Lemma 3 and the fact that the union of two small sets is small ([13], Proposition 5.5.5). The conclusion follows. $\square$

PROOF OF THEOREM 3.2. Since $(X_t)$ is geometrically ergodic, it is $\beta$-mixing, with $E_\pi g(X_t) = E_\pi \varepsilon_t^{2r} < \infty$. It follows that $(\varepsilon_t^2)$ and $(\sigma_t)$ are $\beta$-mixing processes. The fact that $\varepsilon_t$ inherits this $\beta$-mixing property follows from the independence between $\sigma_t$ and $\eta_t$ (see, e.g., [10], proof of Theorem 3). $\square$

3.4. *Existence of moments.* Theorem 3.2 ensures the existence of a moment of some order $2r$. For statistical applications, however, it is often necessary to assume second order stationarity or the existence of higher order moments. The following theorem provides a sufficient condition for the existence of $2p$th-order moments:

THEOREM 3.3. *Let $p \in \mathbb{N}$. Under Assumption A, with $\mu_{2p} = E\eta_t^{2p} < \infty$, if*

$$(3.7) \qquad 0 \leq \alpha < \max_{m \in \{1,2,\ldots\}} \left( \frac{k^{m-1}}{\mu_{2p}\mu_{2m}^{1-1/m}\mu_{2mp}^{1/m}} \right)^{1/(2p+m-1)},$$

*then there exists a strictly stationary solution process $(\varepsilon_t)$ to model (1.2) such that $E_\pi(\varepsilon_t^{2p}) < \infty$.*

REMARK 6. For $m = 1$, the term inside the brackets reduces to $\mu_{2p}^{-1/p}$. A simple condition for the existence of $E(\varepsilon_t^{2p})$ is thus

$$(3.8) \qquad \mu_{2p}\alpha^p < 1,$$

which is also necessary in the standard ARCH(1) case ($k = 0$). However, the example below shows that when $k$ increases, the upper bound in (3.7) is attained for integers $m > 1$.

PROOF OF THEOREM 3.3. Following the same approach as that used in the proof of Lemma 4, but now with $g(x) = g(x_1, x_2) = x_1^p$, we get

$$E[\varepsilon_{t+2}^{2p}|X_t = (x_1, x_2)]$$
$$= E[\eta_{t+2}^{2p}\psi\{\psi(x_1,x_2)\eta_{t+1}^2, x_1\}^p|X_{t-1} = (x_1, x_2)]$$
$$= \mu_{2p}E\{\omega + \alpha\psi(x_1,x_2)\eta_{t+1}^2 \mathbb{1}_{\eta_{t+1}^2 > kx_1/\psi(x_1,x_2)}\}^p$$
$$(3.9) \qquad = \mu_{2p}\sum_{s=0}^{p}\binom{s}{p}\omega^{p-s}\alpha^s\psi(x_1,x_2)^s E[\eta_{t+1}^{2s}\mathbb{1}_{\eta_{t+1}^2 > kx_1/\psi(x_1,x_2)}].$$

By the Hölder and Markov inequalities, we have, for $m \geq 1$,

$$E(\eta_{t+1}^{2s}\mathbb{1}_{\eta_{t+1}^2 > kx_1/\psi(x_1,x_2)}) \leq \{E(\eta_{t+1}^{2ms})\}^{1/m}\left\{P\left[\eta_{t+1}^{2m} > \left(\frac{kx_1}{\psi(x_1,x_2)}\right)^m\right]\right\}^{(m-1)/m}$$
$$\leq \{E(\eta_{t+1}^{2ms})\}^{1/m}\left\{\frac{E(\eta_{t+1}^{2m})\psi(x_1,x_2)^m}{(kx_1)^m}\right\}^{(m-1)/m}$$
$$= \mu_{2ms}^{1/m}\mu_{2m}^{(m-1)/m}\left\{\frac{\psi(x_1,x_2)}{kx_1}\right\}^{m-1}.$$

ERGODICITY OF NONLINEAR ARCH 13When $x_1 > kx_2$ and $x_1 \to +\infty$, the right-hand side of (3.9) is thus bounded by a term which is equivalent to

$$\mu_{2p}\alpha^{2p}\mu_{2mp}^{1/m}\mu_{2m}^{(m-1)/m}\left\{\frac{\alpha}{k}\right\}^{m-1}x_1^p$$

$$=\left\{\alpha\frac{\mu_{2p}^{1/(2p+m-1)}\mu_{2mp}^{1/(m(2p+m-1))}\mu_{2m}^{(m-1)/(m(2p+m-1))}}{k^{(m-1)/(2p+m-1)}}\right\}^{2p+m-1}x_1^p.$$

The right-hand side term inside the brackets being, in view of (3.7), strictly less than 1 for some $m \geq 1$, we thus have

$$\mu_{2p}\alpha^{2p}\mu_{2mp}^{1/m}\mu_{2m}^{(m-1)/m}\left\{\frac{\alpha}{k}\right\}^{m-1}x_1^p \leq (1-\beta)x_1^p - \beta$$

for some constant $\beta > 0$. Therefore, there exists $M > 0$ such that

$$x_1 > M \quad \text{and} \quad x_1 > kx_2 \implies E[\varepsilon_{t+2}^{2p}|X_t = (x_1, x_2)] \leq (1-\beta)x_1^p - \beta.$$

Furthermore, for $x_1 \leq M$, we have $\psi(x_1, x_2) \leq \omega + \alpha M$ and thus, in view of (3.9),

$$E[\varepsilon_{t+2}^{2p}|X_t = (x_1, x_2)] \leq \mu_{2p}\sum_{s=0}^{p}\binom{s}{p}\omega^{p-s}\alpha^s(\omega + \alpha M)^s$$

$$= \mu_{2p}\{\omega + \alpha(\omega + \alpha M)\}^p.$$

Finally, for $x_1 \leq kx_2$, we have $\psi(x_1, x_2) = \omega$ and thus, from (3.9),

$$E[\varepsilon_{t+2}^{2p}|X_t = (x_1, x_2)] \leq \mu_{2p}\{\omega + \alpha\omega\}^p.$$

We can conclude that

$$E[\varepsilon_{t+2}^{2p}|X_t = (x_1, x_2)] \leq \begin{cases} (1-\beta)x_1^{2p} - \beta, & x \in C^c, \\ b, & x \in C, \end{cases}$$

for some strictly positive constants $\beta$ and $b$, with $C = [0, M]^2 \cup \{x_1 \leq kx_2\}$. The theorem follows. $\square$

When $k \leq 1$, a necessary condition can be straightforwardly obtained as follows. Let $(\varepsilon_t)$ be a strictly stationary solution of model (1.2) with a finite $2p$th moment. Then

$$E(\varepsilon_t^{2p}) \geq \mu_{2p}[\omega^p + \alpha^p E\{\varepsilon_{t-1}^{2p}\mathbb{1}_{\varepsilon_{t-1}^2 > k\varepsilon_{t-2}^2}\}]$$

$$= \mu_{2p}[\omega^r + \alpha^p E(\varepsilon_t^{2p}) - \alpha^p E\{\varepsilon_{t-1}^{2p}\mathbb{1}_{\varepsilon_{t-1}^2 \leq k\varepsilon_{t-2}^2}\}]$$

$$\geq \mu_{2p}[\omega^p + \alpha^p E(\varepsilon_t^{2p}) - \alpha^p k^p E(\varepsilon_{t-2}^{2p})].$$

It follows that

$$\{1 - \mu_{2p}\alpha^p(1-k^p)\}E(\varepsilon_t^{2p}) \geq \mu_{2p}\omega^p.$$

TABLE 1
*Constraints for the existence of the second order moment ($p = 1$), for the standard normal distribution, as functions of $k$. The second column gives the value of $m$ for which the maximum is attained in (3.7). The third column gives the constraint for $\alpha$ as a function of $k$ and the last column gives the maximum value for $\alpha$ when $k$ is equal to the upper bound of the interval*

| $k$ | $m$ | $\alpha$ | $\alpha_{\max}$ |
|---|---|---|---|
| $[0, 3[$ | 1 | $[0, 1[$ | 1 |
| $[3, 6.455[$ | 2 | $[0, \{\frac{k}{\mu_4}\}^{1/3}[$ | 1.291 |
| $[6.455, 12.652[$ | 3 | $[0, \{\frac{k^2}{\mu_6}\}^{1/4}[$ | 1.807 |
| $[12.652, 23.714[$ | 4 | $[0, \{\frac{k^3}{\mu_8}\}^{1/5}[$ | 2.635 |
| $[23.714, 43.297[$ | 5 | $[0, \{\frac{k^4}{\mu_{10}}\}^{1/6}[$ | 3.936 |
| $[43.297, 77.694[$ | 6 | $[0, \{\frac{k^5}{\mu_{12}}\}^{1/7}[$ | 5.976 |
| $[77.694, 137.715[$ | 7 | $[0, \{\frac{k^6}{\mu_{14}}\}^{1/8}[$ | 9.181 |

Therefore, a necessary condition for $E(\varepsilon_t^{2p}) < \infty$ is

$$\mu_{2p}\alpha^p(1 - k^p) < 1, \tag{3.10}$$

and we have

$$E(\varepsilon_t^{2p}) \geq \frac{\mu_{2p}\omega^p}{1 - \mu_{2p}\alpha^r(1 - k^p)}.$$

When $k = 0$, (3.10) coincides with (3.8) and provides the necessary and sufficient condition for the existence of $E(\varepsilon_t^{2p})$ in the standard ARCH(1) case (see [12] for moment conditions for the GARCH($p, q$) model).

EXAMPLE. In the case of the standard $\mathcal{N}(0, 1)$ distribution for $\eta_t$, condition (3.7) can be made explicit. First, let $p = 1$. We have $\mu_{2m} = \frac{(2m)!}{2^m m!}$ and simple algebra shows that the maximum in (3.7) is attained for

$$m_0 = m_0(k) = \min_{m \in \{2, 3, \ldots\}} \left\{ m : k < \left( \frac{(2m-1)^m}{\mu_{2(m-1)}} \right)^{1/2} \right\} - 1.$$

Thus, the second-order stationarity condition is

$$0 < \alpha < \left( \frac{k^{m_0 - 1} 2^{m_0} m_0!}{(2m_0)!} \right)^{1/(m_0+1)}.$$

For $k \geq 3$, values of $\alpha$ that are greater than 1 can be compatible with second-order stationarity, as can be seen from Table 1.

Similar computations can be carried out when $p = 2$. Table 2 provides the fourth-order stationarity constraints, for different ranges of values of $k$. The values of $m$ corresponding to the maximum in (3.7) have been obtained numerically. For $k < 3.416$, the maximum is reached for $m = 1$ and the constraint is that of a standard ARCH(1) ($3\alpha^2 < 1$). Interestingly, when $k$ increases, the maximum is reached for larger values of $m$ (e.g., $m = 2$ for $1.763 \leq k < 1.886$) and larger values for $\alpha$ are obtained. It is seen that values of $\alpha$ much larger than 1 are compatible with $E(\varepsilon_t^4) < \infty$ when $k$ is large. Similar tables can be constructed for any value of $p$ and for other distributions.

TABLE 2
*As in Table 1, but for the moment of order 4 ($p = 2$)*

| $k$ | $m$ | $\alpha$ | $\alpha_{\max}$ |
|---|---|---|---|
| $[0, 3.416)$ | 1 | $[0, \frac{1}{\mu_4^{1/2}})$ | 0.577 |
| $[3.416, 4.579)$ | 2 | $[0, \{\frac{k}{\mu_4^{3/2} \mu_8^{1/2}}\}^{1/5})$ | 0.612 |
| $[4.579, 6.373)$ | 3 | $[0, \{\frac{k^2}{\mu_4 \mu_6^{2/3} \mu_{12}^{1/3}}\}^{1/6})$ | 0.684 |
| $[6.373, 8.846)$ | 4 | $[0, \{\frac{k^3}{\mu_4 \mu_8^{3/4} \mu_{16}^{1/4}}\}^{1/7})$ | 0.787 |
| $[8.846, 12.183)$ | 5 | $[0, \{\frac{k^4}{\mu_4 \mu_{10}^{4/5} \mu_{20}^{1/5}}\}^{1/8})$ | 0.923 |
| $[12.183, 16.656)$ | 6 | $[0, \{\frac{k^5}{\mu_4 \mu_{12}^{5/6} \mu_{24}^{1/6}}\}^{1/9})$ | 1.098 |
| $[16.656, 22.626)$ | 7 | $[0, \{\frac{k^6}{\mu_4 \mu_{14}^{6/7} \mu_{28}^{1/7}}\}^{1/10})$ | 1.320 |
| $[22.626, 30.571)$ | 8 | $[0, \{\frac{k^7}{\mu_4 \mu_{16}^{7/8} \mu_{32}^{1/8}}\}^{1/11})$ | 1.599 |
| $[30.571, 41.122)$ | 9 | $[0, \{\frac{k^8}{\mu_4 \mu_{18}^{8/9} \mu_{36}^{1/9}}\}^{1/12})$ | 1.948 |

The outputs of Tables 1 and 2 are represented in Figure 1.

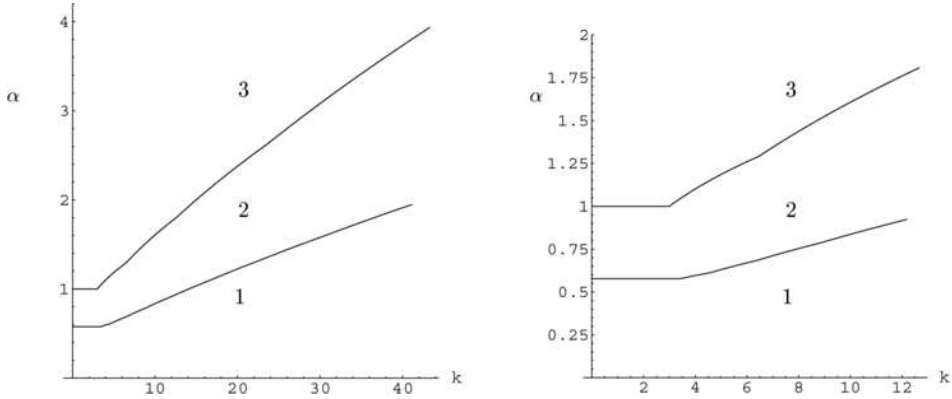

FIG. 1. *Stationarity regions for model* (1.2) *with* $\eta_t \sim \mathcal{N}(0, 1)$. 1. *Existence of* $E\varepsilon_t^4$; 2. *Existence of* $E\varepsilon_t^2$ *with* $E\varepsilon_t^4 = \infty$; 3. *Strict stationarity with* $E\varepsilon_t^2 = \infty$. *The right panel is a zoom of the left panel.*

**Acknowledgment.** The authors wish to thank the referee for careful reading and valuable comments.

Université Mohammed V-Souissi
Angle avenue El Fassi et Cherkaoui
B.P. 8007. N.U. Rabat
Marocco
E-mail: saidiyoussef@hotmail.com

CREST
3 avenue Pierre Larousse
92245 Malakoff Cedex
France
E-mail: zakoian@ensae.fr